\documentclass[12pt]{article}
\usepackage{amsfonts}

\usepackage{mathrsfs}
\usepackage{amscd}
\usepackage{amsmath,amsfonts,amssymb,amscd}
\usepackage{indentfirst,graphics,epsfig,psfrag,color}
\input{epsf}

\usepackage{amsmath,amssymb,amsthm, amscd, amsfonts, graphicx,ifpdf}
\usepackage{lineno}

\newtheorem{thm}{Theorem}[section]

\newtheorem{lem}{Lemma}[section]

\def\qed{\nopagebreak\hfill{\rule{4pt}{7pt}}
\medbreak}

\def\pf{\noindent {\it Proof.} }

\title{\bf\Large Note on Path-Connectivity \\of Complete Bipartite Graphs}
\author{
\small Shasha Li$^1$, Yan~Zhao$^2$\\
\small $^1$ School of Mathematics and Statistics, Ningbo University, Ningbo 315211, Zhejiang, China\\
\small $^2$ Department of Mathematics, Taizhou University, Taizhou 225300, Jiangsu, China\\
\small Email: lishasha@nbu.edu.cn; zhaoyan81.2008@163.com
}
\date{}
\begin{document}

\maketitle

\begin{abstract}
For a graph $G=(V,E)$ and a set $S\subseteq V(G)$ of size at least $2$,
a path in $G$ is said to be an {\it $S$-path} if it connects all vertices of $S$.
Two $S$-paths $P_1$ and $P_2$ are said to be {\it internally
disjoint} if $E(P_1)\cap E(P_2)=\emptyset$ and $V(P_1)\cap V(P_2)=S$.
Let $\pi_G (S)$ denote the maximum number
of internally disjoint $S$-paths in $G$. The {\it $k$-path-connectivity} $\pi_k(G)$
of $G$ is then defined as the minimum $\pi_G (S)$, where $S$ ranges over all $k$-subsets of $V(G)$.
In [M. Hager, Path-connectivity in graphs, Discrete Math. 59(1986), 53--59],
the $k$-path-connectivity of the complete bipartite graph $K_{a,b}$ was
calculated, where $k\geq 2$. But, from his proof,
only the case that $2\leq k\leq min\{a,b\}$ was considered. In this paper, we calculate the
the situation that $min\{a,b\}+1\leq k\leq a+b$ and complete the result.
\\
{\bf Keywords:} path-connectivity; internally disjoint paths; complete bipartite graphs  \\
{\bf AMS Subject Classification 2010:} 05C38, 05C40.
\end{abstract}

\section{Introduction}
The connectivity is one of the most basic concepts of
graph-theoretic subjects, both in a combinatorial sense and an algorithmic sense.
The {\it connectivity} of a graph $G$, denoted by $\kappa(G)$, is defined as the minimum cardinality of
a subset $Q$ of vertices such that $G-Q$ is disconnected or trivial.
A graph $G$ is {\it $k$-connected} if $\kappa(G)\geq k$.
A well-known result by Whitney \cite{Whitney} provided an equivalent definition of
connectivity as follows: A graph $G$ is $k$-connected if and only if there exist $k$ internally disjoint paths
between any two vertices in $G$.

Dirac \cite{Dirac} showed that
in a $(k-1)$-connected graph, there is a path through each $k$ vertices.
In order to generalize this conclusion, the concept of path-connectivity was introduced
by Hager \cite{Hager2}. Given a graph $G=(V,E)$ and a set $S\subseteq V(G)$ of size at least $2$,
a path in $G$ is said to be an {\it $S$-path} if it connects all vertices of $S$.
Two $S$-paths $P_1$ and $P_2$ are said to be {\it internally
disjoint} if $E(P_1)\cap E(P_2)=\emptyset$ and $V(P_1)\cap V(P_2)=S$.
Let $\pi_G (S)$ denote the maximum number
of internally disjoint $S$-paths in $G$. The {\it $k$-path-connectivity}
of $G$, denoted by $\pi_k(G)$, is then defined
as $\pi_k(G)=$min$\{\pi_G(S)| S\subseteq V(G)\ and\ |S|=k\}$,
where $2\leq k\leq n$. Clearly, $\pi_2(G)$ is
exactly the classical connectivity $\kappa(G)$ and $\pi_n(G)=\pi_G(V)$ is exactly
the maximum number of edge-disjoint spanning paths in $G$. Moreover,
by the concept of path-connectivity, the conclusion of Dirac can be restated as:
if $\kappa(G)\geq k-1$, then $\pi_k(G)\geq 1$.

In \cite{Hager2}, Hager studied the sufficient conditions for $\pi_k(G)$
to be at least $\ell$ in terms of $\kappa(G)$. In addition, he calculated the exact values of
$\pi_k(G)$ for complete graphs and complete bipartite graphs.

\begin{lem}[\cite{Hager2}]
 Let $G$ be a complete graph of order $n$. Then
\begin{center}
$\pi_k(K_n)=\lfloor\frac{2n+k^2-3k}{2(k-1)}\rfloor$, for $2\leq k\leq n$.
\end{center}
\end{lem}

\begin{lem}[\cite{Hager2}] \label{lem6}
Let $K_{a,b}$ be a complete bipartite graph with parts of sizes $a$ and $b$. Then
\begin{center}
$\pi_k(K_{a,b})=min\{\lfloor\frac{a}{k-1}\rfloor,\lfloor\frac{b}{k-1}\rfloor\}$, for $k\geq2$.
\end{center}
\end{lem}

However, the result for complete bipartite graphs is incomplete. Actually, from the proof of Lemma \ref{lem6} in \cite{Hager2},
only the case that $2\leq k\leq min\{a,b\}$ was considered. In this paper,
we calculate the missing cases and obtain the following result:

\begin{thm}\label{thm1}
Let $K_{a,b}$ be a complete bipartite graph with parts of sizes $a$ and $b$, where $1\leq a\leq b$. \\
If $2\leq k\leq a$, then $\pi_k(K_{a,b})=\lfloor\frac{a}{k-1}\rfloor$;\\
otherwise, if $a=b$, then $\pi_k(K_{a,b})=
\begin{cases}
1,\ a=3\ and\ k=4,\ or\ a+1\leq k\leq 2a-3(a\geq 4),\\
2,\ k=2a-2(a\geq 4),\\
max\{\lfloor\frac{a}{2}\rfloor,1\},\ k=2a-1(a\geq 2)\ or\ k=2a(a\geq 1);
\end{cases}
$

$\qquad \quad $if $b=a+1(a\geq 1)$, then $\pi_k(K_{a,b})=
\begin{cases}
1,\ a+1\leq k\leq 2a,\\
\lfloor\frac{a+1}{2}\rfloor,\ k=2a+1=a+b;
\end{cases}
$\\

$\qquad \quad $if $b\geq a+2(a\geq 1)$, then $\pi_k(K_{a,b})=
\begin{cases}
1,\ k=a+1,\\
0,\ a+2\leq k\leq a+b.
\end{cases}
$
\end{thm}

\subsection{Preliminaries}

Before proving the theorem, we first introduce some concepts.
All graphs considered here are finite and simple.
For a vertex $v$ in a graph $G$, denote by $N_G(v)$ the set of neighbors of $v$ in $G$. For a subset
$U \subseteq V(G)$, the subgraph induced by $U$ is denoted by $G[U]$.
We write $N_U(v)$ instead of $N_{G[U]}(v)$.
A $spanning$ $subgraph$ of a graph $G$ is a subgraph that contains all the vertices of $G$.
We refer the reader to \cite{Bondy} for the terminology and notations not defined in this paper.

Next, we present the maximum number of edge-disjoint
spanning paths in a complete bipartite graph that will be used in the proof of Theorem \ref{thm1}.

\begin{lem}\label{lem7}
Let $K_{a,b}$ be a complete bipartite graph with parts of sizes $a$ and $b$, where $1\leq a\leq b$.
If $1\leq a\leq b\leq a+1$, then the maximum number $\pi_{a+b}(K_{a,b})$ of edge-disjoint
spanning paths in $K_{a,b}$ is $max\{\lfloor\frac{b}{2}\rfloor,1\}$; otherwise, $\pi_{a+b}(K_{a,b})=0$.
\end{lem}
\pf Obviously, if $b\geq a+2$, there is no spanning path in $K_{a,b}$. Next, we
consider the following two cases: $b=a$ and $b=a+1$, respectively. Suppose that
$X$ and $Y$ are two parts of $K_{a,b}$, where $X=\{x_1,x_2,\ldots,x_a\}$ and
$Y=\{y_1,y_2,\ldots,y_b\}$.

{\it Case 1:} $b=a+1$.

Let $P_1$ be a path such that $N_{P_1}(x_i)=\{y_i,y_{i+1}\}$, where $1\leq i\leq a$,
i.e. $P_1=y_1x_1y_2x_2\ldots x_ay_{a+1}$; let $P_2$ be a path such that
$N_{P_2}(x_i)=\{y_{i+2},y_{i+3}\}$, where $1\leq i\leq a$ and the subscripts of $y$
are taken module $b$, i.e. $P_2=y_3x_1y_4x_2\ldots x_ay_2$, and so on.
Therefore, in this way, we can find at least $\lfloor\frac{b}{2}\rfloor$
edge-disjoint spanning paths $P_j$ with $N_{P_j}(x_i)=\{y_{i+2(j-1)},y_{i+2j-1}\}$,
where $1\leq j\leq \lfloor\frac{b}{2}\rfloor$, $1\leq i\leq a$ and the subscripts of $y$
are taken module $b$.

On the other hand, the number of edge-disjoint spanning paths in  $K_{a,b}$ is
at most $\lfloor\frac{ab}{a+b-1}\rfloor=\lfloor\frac{(b-1)b}{2b-2}\rfloor=\lfloor\frac{b}{2}\rfloor$,
since each spanning path needs $a+b-1$ edges.

Hence, in this case, $\pi_{a+b}(K_{a,b})=\lfloor\frac{b}{2}\rfloor$.

{\it Case 2:} $b=a$.

If $a=b=1$, obviously $\pi_2(K_{1,1})=1$.

For $a=b\geq 2$, let $P_1$ be a path such that $N_{P_1}(x_i)=\{y_i,y_{i+1}\}$($1\leq i\leq a-1$)
and $N_{P_1}(x_a)=\{y_a\}$; let $P_2$ be a path such that
$N_{P_2}(x_i)=\{y_{i+2},y_{i+3}\}$($1\leq i\leq a-1$ and the subscripts
are taken module $a$) and $N_{P_2}(x_a)=\{y_2\}$, and so on. Therefore, in this way, we can
obtain at least $\lfloor\frac{a}{2}\rfloor=\lfloor\frac{b}{2}\rfloor$
edge-disjoint spanning paths $P_j$ with $N_{P_j}(x_i)=\{y_{i+2(j-1)},y_{i+2j-1}\}$ and
$N_{P_j}(x_a)=\{y_{a+2j-2}\}$, where $1\leq j\leq \lfloor\frac{b}{2}\rfloor$, $1\leq i\leq a-1$
and the subscripts are taken module $a$.

On the other hand, since $2a^2< 2a^2+a-1=(2a-1)(a+1)$, $\frac{a^2}{2a-1}< \frac{a+1}{2}$
and so $\lfloor\frac{a^2}{2a-1}\rfloor\leq \lfloor\frac{a}{2}\rfloor$. Thus, we have at most
$\lfloor\frac{a}{2}\rfloor=\lfloor\frac{b}{2}\rfloor$ edge-disjoint spanning paths,
i.e. $\pi_{a+b}(K_{a,b})=\lfloor\frac{b}{2}\rfloor$. The proof is complete.\qed

\section{Proof of Theorem \ref{thm1}.}

Let $K_{a,b}=G$ and suppose that
$X$ and $Y$ are two parts of $K_{a,b}$, where $X=\{x_1,x_2,\ldots,x_a\}$ and
$Y=\{y_1,y_2,\ldots,y_b\}$.

For $2\leq k\leq a$, the result was given by Lemma \ref{lem6}. Thus, we consider only
the situation that $a+1\leq k\leq a+b$. We distinguish three cases as follows.

{\it Case 1:} $a=b$.

If $k=a+b=2a$, by Lemma \ref{lem7}, $\pi_k(K_{a,b})=max\{\lfloor\frac{a}{2}\rfloor,1\}$
for all $a\geq 1$.

For $k=2a-1(a\geq 2)$, let $S$ consist of any $k$ vertices of $V(G)$. Since
$K_{a,b}[S]=K_{a,a-1}$, $\pi_G(S)\geq \pi_{K_{a,a-1}}(S)=\pi_{2a-1}(K_{a,a-1})=
max\{\lfloor\frac{a}{2}\rfloor,1\}=\lfloor\frac{a}{2}\rfloor$ by Lemma \ref{lem7}. Therefore,
in $K_{a,b}[S]$, we have $\lfloor\frac{a}{2}\rfloor$ internally disjoint $S$-paths.
Now, there are $a(a-1)-\lfloor\frac{a}{2}\rfloor(2a-2)$
edges left in $K_{a,b}[S]$ and a vertex left in $V(G)\setminus S$.
Note that, if there is only one vertex available in $V(G)\setminus S$,
we need $k-2=2a-3$ edges of $K_{a,b}[S]$ more to form an $S$-path.
But, since $a(a-1)-\lfloor\frac{a}{2}\rfloor(2a-2)$ is equal to $0$ if $a$ is even and
is equal to $a-1$ if $a\geq 3$ is odd, the remaining vertex and edges are not enough to form an
$S$-path. Thus, $\pi_G(S)=\lfloor\frac{a}{2}\rfloor$. By the arbitrariness of $S$,
$\pi_k(K_{a,b})=\lfloor\frac{a}{2}\rfloor$.

Let now $k=2a-2$. If $a=b=3$, then $K_{a,b}=K_{3,3}$ and $k=4$.
Let $S=\{x_1,x_2,x_3,y_1\}$. It is easy to check that $\pi_{K_{3,3}}(S)\leq 1$.
On the other hand, there exists a spanning path in $K_{3,3}$, which is also a path connecting
any four vertices of $V(G)$. Hence, $\pi_4(K_{3,3})=1$.

If $a=b\geq 4$ and $k=2a-2$, let $S$ consist of $k$ vertices of $V(G)$. If
$|S\cap X|=|S\cap Y|=a-1$, then $K_{a,b}[S]=K_{a-1,a-1}$. So by Lemma \ref{lem7}
$\pi_G(S)\geq \pi_{2a-2}(K_{a-1,a-1})=max\{\lfloor\frac{a-1}{2}\rfloor,1\}\geq 2$ for $a\geq 5$
and for $a=4$, it is easy to check that $\pi_G(S)\geq 2$. If
$|S\cap X|=a$ and $|S\cap Y|=a-2$ or $|S\cap X|=a-2$ and $|S\cap Y|=a$,
without loss of generality, let $S=\{x_1,x_2,\ldots,x_a,y_1,y_2,\ldots,y_{a-2}\}$
and $\hat{S}=S\cup \{y_{a-1}\}$. Since $K_{a,b}[\hat{S}]=K_{a,a-1}$ and
by Lemma \ref{lem7} $\pi_{2a-1}(K_{a,a-1})=max\{\lfloor\frac{a}{2}\rfloor,1\}\geq 2$
for $a\geq 4$, we can suppose that $P_1$ and $P_2$ are two internally disjoint $\hat{S}$-paths
in $K_{a,b}[\hat{S}]$. Now, for $P_2$, replace the vertex $y_{a-1}$ by $y_a$ and
replace the two edges incident to $y_{a-1}$ by edges incident to $y_a$. Denote the resulting path by $\hat{P_2}$.
Obviously, $P_1$ and $\hat{P_2}$ are two internally disjoint $S$-paths
in $K_{a,b}$ and so $\pi_G(S)\geq 2$. On the other hand, any path through the $a$
vertices of $X$ needs $a-1$ vertices of $Y$ and so any $S$-path needs a vertex from
$V(G)\backslash S=\{y_{a-1},y_a\}$. Hence, in this case $\pi_G(S)=2$.
It follows that $\pi_{2a-2}(K_{a,a})=2$ for $a\geq 4$.

Let now $a+1\leq k\leq 2a-3$. Since $a+1\leq 2a-3$, we have $a\geq 4$. For $K_{a,a}$, clearly there
is a spanning path, which is also a path connecting any $k$ vertices of $V(G)$.
Now let $S=\{x_1,x_2,\ldots,x_a,y_1,\ldots,y_{k-a}\}$. Since $S\supset X$, each
$S$-path needs at least $a-1$ vertices of $Y$ and so needs
$a-1-(k-a)=2a-k-1$ vertices from $V(G)\backslash S=\{y_{k-a+1},\ldots,y_a\}$.
Since $|V(G)\backslash S|=2a-k$ and $2a-k-1\geq 2a-(2a-3)-1=2$,
there exists at most one $S$-path. Hence, $\pi_k(K_{a,b})=1$
in this case.

{\it Case 2:} $b=a+1$.

For $a+1\leq k\leq a+b-1=2a(a\geq 1)$, clearly $\pi_k\geq 1$.
Let $S$ be a subset of $V(G)$ with $|S|=k$ and $S\supseteq Y$.
Then any $S$-path must contain all the vertices of $X$ and so
must be a spanning path of $G$. Since $V(G)\backslash S\neq \emptyset$ and
the vertices in $V(G)\backslash S$ can be used only once, $\pi_G(S)\leq 1$
and so $\pi_k(K_{a,b})=1$.

If $k=a+b=2a+1$, $\pi_k(K_{a,b})$ is exactly the maximum number of
edge-disjoint spanning paths in $K_{a,b}$. So by Lemma \ref{lem7},
$\pi_k(K_{a,b})=max\{\lfloor\frac{b}{2}\rfloor,1\}=\lfloor\frac{a+1}{2}\rfloor$
for all $a\geq 1$.

{\it Case 3:} $b\geq a+2$.

If $k=a+1(a\geq 1)$, for any $k$-subset $S$ of $V(G)$, there exists
a subset $\hat{S}$ of $V(G)$ such that $\hat{S}\supset S$, $\hat{S}\supset X$
and $|\hat{S}\cap Y|=a+1$. Clearly, there is a path connecting $\hat{S}$ and
also connecting $S$. So, $\pi_k(K_{a,b})\geq 1$. On the other hand, if $|S\cap Y|=k$,
i.e. $S\cap X=\emptyset$, then obviously $\pi_G(S)\leq 1$. Hence,
$\pi_k(K_{a,b})=1$ in this case.

For $a+2\leq k\leq a+b$, let $S$ be a subset of $V(G)$ such that $|S|=k$
and $|S\cap Y|\geq a+2$. Since $|X|=a$, $\pi_G(S)=0$ and so $\pi_k(K_{a,b})=0$ in this case.

We have considered all cases and so the proof is complete.\qed

\noindent{\bf Acknowledgments.} Shasha Li was
supported by Zhejiang Provincial Natural Science Foundation of China (No. LY18A010002)
and the Natural Science Foundation of Ningbo, China.
Yan Zhao was supported by the National Natural Science Foundation of China (No. 11901426) and Qing Lan Project.

\end{document}